\documentclass[10pt,leqno]{amsart}
\usepackage{amssymb}
\usepackage{mathrsfs}
\usepackage{color}
\usepackage{soul}

\newcommand\R{{\mathbb R}}

\newcommand\EEE{{\bf E}}
\newcommand\PPP{{\bf P}}

\def\BB{{\mathcal B}}

\def\FF{{\mathcal F}}

\def\KK{{\mathcal K}}
\def\LL{{\mathcal L}}

\def\NN{{\mathcal N}}
\def\OO{{\mathcal O}}

\def\SS{{\mathcal S}}

\def\VV{{\mathcal V}}
\def\WW{{\mathcal W}}

\def\YY{{\mathcal Y}}
\def\ZZ{{\mathcal Z}}

\def\PPP{{\mathbf P}}

\def\Ps{{\bf P}_{\! sym}}

\def\SIGMA{\mathfrak{S}}
\def\BSN{{\mathcal {BS}^N}}

\def\eps{{\varepsilon}}

\newcommand{\wto}{\rightharpoonup}

\newtheorem{theo}{Theorem}

\newcommand{\beqn}{\begin{equation}}
\newcommand{\eeqn}{\end{equation}}
\newcommand{\bear}{\begin{eqnarray}}
\newcommand{\eear}{\end{eqnarray}}
\newcommand{\bean}{\begin{eqnarray*}}
\newcommand{\eean}{\end{eqnarray*}}

\newcommand{\Black}{\color{black}}

\def\signsm{\bigskip \begin{center} {\sc St\'ephane Mischler\par\vspace{3mm}
Universit\'e Paris-Dauphine \& IUF \par
CEREMADE, UMR CNRS 7534\par
Place du Mar\'echal de Lattre de Tassigny
75775 Paris Cedex 16\par
FRANCE\par\vspace{3mm}
e-mail:} \tt{mischler@ceremade.dauphine.fr} \end{center}}

\begin{document}
 
 \title[Kac's chaos and Kac's program]
{Kac's chaos and Kac's program \\ 
}

\author{S. Mischler}

\begin{abstract} In this note I present the main results about the quantitative and qualitative propagation of chaos for the Boltzmann-Kac system obtained  in collaboration with C.~Mouhot in \cite{MMinvent} which  gives a possible answer to some questions formulated by Kac in \cite{Kac1956}.  We also present some related recent results about Kac's chaos and Kac's program obtained in \cite{MMWchaos,HaurayMischler,KleberSphere} by  K.~Carrapatoso, M. Hauray, C. Mouhot, B. Wennberg 
and myself.

\end{abstract}

\maketitle

\begin{center} {\bf Manuscript version of a talk given in 
}

\smallskip
{\bf   s\'eminaire Laurent Schwartz, 2012-2013  
}

\end{center}


Jump processes.

%
%
%
%
%
	
\vspace{0.3cm}



\bigskip


\tableofcontents



\section{Introduction} 
\label{sec:intro}
\setcounter{equation}{0}
\setcounter{theo}{0}

\subsection
{6th Hilbert Problem}

The Boltzmann equation was introduced by  Maxwell (1867, \cite{Maxwell1867}) 
and Boltzmann (1872, \cite{Boltzmann1872}) in order to describe the evolution of 
a rarefied gas in which particles uniquely interact through binary collisions. 
That equation governs the time evolution of the statistical distribution $f(t,x,v) \ge 0$  of positions and velocities of particles of the gas. 
The Boltzmann equation, together with  the Vlasov equation for collisionless gas,  is the most fundamental model in classical kinetic theory of gases.

\smallskip
In 1900 during the  conference of the International Congress of Mathematicians in Paris, D.~Hilbert  invited us  in his 6th problem 
to develop an axiomatic approach of physics, and in particular of the just beginning kinetic theory of gases, in the following words:  

{\it 
``The investigations on the  {foundations} of geometry suggest the 
problem: To treat in the same manner,  {by means of axioms}, 
those physical sciences in which mathematics plays an 
important part; in the first rank are the  {theory of probabilities 
and mechanics}.
As to the axioms of the theory of probabilities, it seems to me 
desirable that their logical investigation should be accompanied 
by a rigorous and satisfactory development of the method of 
mean values in mathematical physics, and in particular in  {the 
kinetic theory of gases"}.
} 

\smallskip
In other words, for the Boltzmann equation
 Hilbert's question is the following: is-it possible to obtain the Boltzmann equation at the statistical level of description from a microscopic description of the gas dynamic, namely from the dynamics of molecules  governed by the Newton's law of motions? 

The ``Boltzmann-Grad limit" which explain how to get the Boltzmann equation from such a microscopic description 
was identified by Grad (1958, \cite{Grad1958}) ({thanks to the BBGKY method}) and mathematically rigorously proved by Lanford (1975, \cite{Lanford}) on a very small time interval (smaller that the necessary waiting time  before half of all the particles collide once). Of course the ``Boltzmann-Grad limit" is very interesting and very difficult to justify, in particular because it requires to understand how to get an irreversible equation (the Boltzmann equation) from a reversible equation (the Newton's law of motions), and very few results are known on that major problem up to now.
We refer to \cite{BGLS} and the references therein for updated results on that direction.

\subsection{Kac's approach}

In order to  circumvent the above difficulties, M. Kac (1956, \cite{Kac1956}) suggested to derive the space homogeneous Boltzmann equation as the limit as the number of particles $N$ goes to infinity of 
a system of $N$ indistinguishable particles which velocities are modified through stochastic  collisions. 

\smallskip
The underlying limit procedure is a mean field limit in the sense that each particle interacts with all the other particles with an intensity of order  
$\OO(1/N)$. That limit  is different from the ``Boltzmann-Grad limit". 

\smallskip
Although that problem is clearly simpler than the justification of the Boltzmann-Grad limit, it is still an interesting and difficult mathematical problem.
It furthermore shares with the ``Boltzmann-Grad limit" the same mathematical difficulty to work with a family of functional spaces with increasing dimension and the necessity to rigorously define the mathematical notion of  asymptotic stochastic independence (Kac's chaos) and
thus to clarify the notion of molecular chaos on which Boltzmann's work is based.

\subsection{$N$-particle system}
The approach by Kac was generalized by McKean (1967, \cite{McKean1967}), and then many other people. 
Generally and roughly speaking, the approach is as follows. 
Consider a system of $N$ indistinguishable particles, each particle being identified by its state (position, velocity)  $\ZZ^N = (\ZZ_1, ..., \ZZ_N)$, $\ZZ_i \in E$, $E= \R^d$,  which evolves accordingly to one of the following equations 
$$
d\ZZ_i = {1 \over N} \sum_{j=1}^N a(\ZZ_i - \ZZ_j) \, dt \qquad \hbox{(ODE)} , 
$$
$$
d\ZZ_i= {1 \over N} \sum_{j=1}^N a(\ZZ_i - \ZZ_j) \, dt + \sqrt{2\nu} d\BB_i \qquad \hbox{(Brownian SDE)} , 
$$
$$
d\ZZ = {1 \over N} \sum_{i,j=1}^N\int_{S^{d-1}}(\ZZ'_{ij} - \ZZ) \,  \tilde b_{ij} \,  d\NN(d\sigma,i,j)  \qquad \hbox{(Boltzmann-Kac/Poisson SDE)}.
$$
Here $a$ is an interaction force field,  $(\BB_i)$ is a family of independent  Brownian motions, $\NN$ is a Poisson measure,
$\ZZ'_{ij} = (\ZZ_1, ..., \ZZ'_i, ... ,\ZZ'_j, ..., \ZZ_N)$ is the system of post-collision velocities after the collision of the velocities pair  $(\ZZ_i,\ZZ_j)$,
$\tilde b_{ij} := \tilde b(z_i-z_j,  \sigma)$ is the collision cross-section. We do not explain more the notations in that first introduction and we will make more precise the definition of the 
Boltzmann-Kac system in section~\ref{subsec:BKsystem} below.

\medskip
The law $G^N(t) := \LL(\ZZ^N)$ of $\ZZ^N$ then satisfies the Master equation (Liouville equation in the deterministic case, backward Kolmogorov equation in the stochastic case)  
$$
\partial_t \langle G^N , \varphi \rangle = \langle G^N ,  \Lambda^N  \varphi \rangle \qquad \forall \, \varphi \in C_b(E^N), 
$$
where the generator  $\Lambda^N$ writes

$$
(\Lambda^N \varphi) (Z) :=  {1 \over N} \sum_{i,j=1}^N a(z_i-z_j) \cdot \nabla_i \varphi   \qquad \hbox{(ODE)} ,
$$
$$
(\Lambda^N \varphi) (Z) :=   {1 \over N} \sum_{i,j=1}^N a(z_i-z_j) \cdot \nabla_i \varphi +  \nu \sum_{i=1}^N \Delta_i \varphi  
 \qquad \hbox{ (SDE)} ,
$$ 
$$
(\Lambda^N \varphi) (Z) = {1 \over N} \sum_{1\le i < j\le N}^N 
  \int_{\mathbb{S}^{d-1}}  \left[\varphi(Z'_{ij}) -
    \varphi(Z)\right] \, \tilde b_{ij} \, {\rm d}\sigma  \qquad \hbox{ (Boltzmann-Kac)}.
$$

\subsection{Nonlinear PDE limit}  The question we are interested in is whether we can identify the possible limit as $N$ tends to $\infty$ of the law $\LL(\ZZ^N_1)$ of one typical particle. More precisely, we aim to prove that  $\LL(\ZZ^N_1) \to f = f(t,dz)$ as $N \to \infty$ where  $f\in C([0,\infty);P(E))$ is the unique solution to the evolution PDE  
$$
\partial_t f  = \hbox{div}_z [ (a* f) f]  \qquad (\hbox{Vlasov}),
$$
$$
\partial_t f  = \hbox{div}_z [ (a* f) f]  + \nu \Delta f \qquad (\hbox{McKean-Vlasov}),
$$
$$
\partial_t f  = \int_{\R^d\times S^{d-1}}  [ f(z')f(v') - f(z) f(v)] \, b \, dz d\sigma \quad \ (\hbox{homogeneous Boltzmann}) ,
$$
depending of the dynamic governing the $N$-particle system. 

\medskip
It is not difficult to figure out (at least formally) why the above PDEs are the correct ones. Indeed, assuming that  
$$
\LL(\ZZ^N_1) \to f = f(t,dz), \quad 
\LL(\ZZ^N_1,\ZZ^N_2) \to g = g(t,dz,dv), 
$$
we may easily pass to the limit in the expression $\langle \Lambda^N,\varphi \rangle$ for a given fixed function $\varphi (Z) = \varphi (z_1)$, 
$\varphi \in C^2_b(E)$. Coming back to   the Master equation, we get  
$$
\partial_t f  =  \hbox{div}_z \Bigl[  \int a(z-v) g(dz,dv)  \Bigr]   \qquad (\hbox{Vlasov}),
$$
$$
\partial_t f  =  \hbox{div}_z \Bigl[  \int a(z-v) g(dz,dv)  \Bigr] + \nu \Delta f  f \qquad (\hbox{McKean-Vlasov}),
$$
$$
\partial_t f  = \int_{\R^d\times S^{d-1}}  [ g(z',v') - g(z,v)] \, \tilde b \, dzd\sigma \quad \ (\hbox{homogeneous Boltzmann}). 
$$
We immediately obtain the Vlasov equation, the McKean-Vlasov equation  and the Boltzmann if we make the additional (molecular chaos) independence assumption  
{ $g (v,z) = f(v) \, f(z)$}. 

\medskip
The above picture is not that easy because for $N$ fixed particles the states $ \ZZ_1(t)$, ..., $\ZZ_N (t)$ are {\bf never} independent for positive time $t > 0$  even if the initial states  $\ZZ_1(0), ..., \ZZ_N(0)$ are assumed to be independent : that is an inherent consequence of the fact that particles do interact!

However, the nonlinear PDE can be obtained as a {\it ``law of large numbers" } for a not independent array of variables in the mean-field limit. That is more demanding than an usual law of large numbers for independent and identically distributed sequence of variables. In order to justify that limit,  we  have to prove the propagation of chaos in the sense that  
$$
\LL (\ZZ^N_1(0),\ZZ^N_2(0)) \to f_0 \otimes f_0 \ \quad \Longrightarrow \ \quad
\LL (\ZZ^N_1(t),\ZZ^N_2(t)) \to f_t \otimes f_t,
$$
where $f_t$ is the solution to the associated nonlinear PDE (formal mean-field limit). 
The above convergence can be proved for all the above mentioned type of model but it always requires to investigate more that the sole law $\LL(\ZZ^N_1)$ of a typical particle. Roughly speaking it requires to study the evolution of a least the law $\LL(\ZZ^N_1,\ZZ^N_2)$ of a typical pair of particles (that was the quantity considered by M. Kac in the seminal article \cite{Kac1956}) or more generally to study the laws  $\LL(\ZZ^N_1,...,\ZZ^N_j)$ as $N\to \infty$ for any fixed $j \ge 1$ (BBGKY method) or really all the particles (as in the coupling method or the empirical measures method).

\subsection{Definition of Kac's chaos } For a infinite system $(\ZZ_i)_{i \ge 1}$ of particles the (molecular) Boltzmann's chaos means that 
$$
\LL(\ZZ_i,\ZZ_j) = f \otimes f, \quad \forall \, i \not= j.
$$
That is the stochastic independence (for a sequence of stochastic variables). 

 \smallskip
 For a sequence $\ZZ^N = (\ZZ^N_1, ..., \ZZ^N_N)$ of $N$-indistinguishable particle systems with $N\to\infty$,  the stochastic chaos according to Kac means that 
$$
\LL (\ZZ^N_i,\ZZ^N_j) \to f \otimes f   \quad\hbox{as}\quad N \to \infty, \quad \forall \, i \not= j \hbox{ fixed}.
$$
That is a kind of asymptotic stochastic independence (of the coordinates of a sequence of stochastic arrays).


\section{Kac's program}
\label{sec:KacProgram}
\setcounter{equation}{0}
\setcounter{theo}{0}

\subsection{Contributions by Kac and Kac's program  }
\  In \cite{Kac1956}, Kac considers a toy model: a $1d$ caricature of the nonlinear Boltzmann equation which is called after him  as ``the Kac's model"
and the associated ``Kac's $N$-particle system".
He defines the mathematical notion of stochastic chaos (for a sequence of stochastic arrays). He proves the propagation of chaos for ``the Kac model" and therefore, for the very first time, he obtain a statistical description of a gas (the "Kac's equation") by passing to the limit  in a miscropic description of that same gas (by mean of the a  ``Kac's $N$-particle system").   He also shows that the $N$-particle dynamic preserves the (so-called) Kac's sphere  
$$
\KK\SS^N := \{ V \in \R^{N}; \, |v_1|^2+   ... + |v_N|^2 = N \},
$$
and that, for any fixed  $N \ge 2$, the $N$-particle system converges in the long time asymptotic  to its equilibrium / invariant measure, namely 
$$
G^N_t = \LL(\VV^N_{1t}, ..., \VV^N_{Nt}) \ \mathop{\longrightarrow}_{t \to \infty} \gamma^N = \hbox{ the uniform measure on  } \KK\SS^N.
$$

\medskip\medskip 
Kac also formulates (in a more or less explicitly manner) a series of questions that we list below and that we refer as the ``Kac's program":  

\smallskip\noindent
{\bf Problem 1.}  Prove the propagation of chaos for some/any  {\it ``realistic"} models; 

\smallskip\noindent
{\bf Problem 2.}  Prove the long time convergence of the $N$-particle system to its equilibrium as $t\to \infty$ with a speed which is uniform with respect 
to the number $N$ of particles;

\smallskip\noindent
{\bf Problem 3.} Establish the $H$-Theorem of Boltzmann for the nonlinear Boltzmann equation 
directly from the microscopic description of  the gas. 
That  last problem seems to be the initial motivation of Kac.

\subsection{The Boltzmann-Kac system and the Boltzmann equation  }
\label{subsec:BKsystem}

Another way to describe the Boltzmann-Kac system from a stochastic trajectories point of view is the following. 
We consider a system of $N$ particles  $\VV^N = (\VV_1, ..., \VV_N)$, $\VV_i \in E= \R^3$, which changes because of stochastic (collisional) jumps, i.e.  
$(\VV^N_t)_{t \ge 0}$ is the Markov process defined (by repeating the process) by   :

\medskip
(i)~for any $(\VV_{i^*},\VV_{j^*})$~one draws a collision time ~$T_{i^*,j^*}~\sim~Exp (B(|\VV_{i^*}- \VV_{j^*}|))$;
and one chooses the (pre-collisional) couple of velocities $(\VV_i,\VV_j)$ such that  
$$
T_{i,j} = \min_{(i^*,j^*)} T_{i^*,j^*} .
$$

\medskip
(ii) one draws an angle  $\sigma \in S^2$ according to the law $b(\cos \, \theta)$,  $\cos\theta = \sigma \cdot u_{ij}$, $u_{ij} = (\VV_i-\VV_j)/|\VV_i-\VV_j|$, and one then defines the post-collisional velocities  $(\VV'_i,\VV'_j)$ by 
\[
\VV'_i = {\VV_i + \VV_j \over 2} + {|\VV_j-\VV_i| \over 2} \, \sigma, \qquad 
\VV'_j = {\VV_i + \VV_j \over 2} - {|\VV_j-\VV_i| \over 2} \, \sigma.
\]

Doing that, the momentum and the (kinetic) energy are conserved during the binary collision  
\[
\VV'_i + \VV'_j = \VV_i + \VV_j, \qquad 
|\VV'_i|^2 + |\VV'_j|^2  = |\VV_i|^2 + |\VV_j|^2,
\]
and therefore also the mean momentum and the mean energy for the system are conserved 
$$
\sum_{i} \VV_i(t) = \hbox{cst}, \quad \sum_{i} |\VV_i(t)|^2 = \hbox{cst}. 
$$
After a change of time scaling, the law $G^N_t  \in \PPP(E^N)$ of the $N$-particle system  $(\VV^N_t)_{t \ge 0}$ satisfies the backward
Kolmogorov  equation  
\beqn\label{eq:BKs}
 \partial_t \langle G^N,\varphi \rangle = \langle G^N, \Lambda^N \varphi \rangle 
\qquad \forall \, \varphi \in C_b(E^N), \quad G^N(0) = G^N_0 , \qquad
\eeqn
where the generator is still given by 
$$
(\Lambda^N\varphi) (V) = {1 \over N} \sum_{i,j= 1}^N  B(v_i-v_j)  \int_{S^2}b(\cos \theta_{ij}) \,  [ \varphi'_{ij} - \varphi] \, d\sigma.
$$
Here we use the shorthand $\varphi = \varphi(V)$,  $\varphi'_{ij} = \varphi(V'_{ij})$, $V'_{ij} = (v_1, .., v'_i, .., v'_j, .., v_N)$. 

\smallskip
In the sequel we only consider the three following classical examples of collisions cross section $\tilde b = B \, b$: 

\begin{itemize}
\item Maxwell interaction with Grad's cutoff  {\bf (MG)}: $B = 1$, $b = 1$;
\item True  Maxwell interaction {\bf (M)}: $B = 1$, $b \notin L^1$;
\item Hard spheres interaction {\bf (HS)}: $B(z) = |z|$, $b=1$. 
\end{itemize} 

The associated nonlinear space homogeneous Boltzmann equation is defined on  $\PPP_2(\R^3)$, 
the space of probability measures with finite second moment,  by
\beqn\label{eq:Beq}
  \partial_t f = Q(f), \quad f(0) = f_0  \qquad
\eeqn
where
$$
\langle Q(f), \varphi \rangle := \int_{\R^6\times S^2} B(v-v_*)  \,  b(\cos\theta)  \, (\phi
(v') -  \phi (v))  \, d\sigma \, f(dv) \, f(dv_*), 
$$
and as before
\[
v' = {v+ v_* \over 2} + {|v-v_*| \over 2} \, \sigma.
\]
The Boltzmann equation generates a nonlinear semigroup
\beqn\label{eq:defSNLt}
\forall \, f_0 \in \PPP_2(\R^3) \qquad S_t^{NL} f_0 := f_t.
\eeqn

\subsection{Problem 1. Known results on the propagation of chaos for the Boltzmann-Kac system}
 
We define the Boltzmann's spheres on $E := \R^3$ by 
$$
\BB\SS^N := \{ V \in E^{N}; \,v_1 +... + v_N = 0, \, |v_1|^2+   ... + |v_N|^2 = N \}.
$$

\begin{theo}\label{theo:Kac1}
Consider  $G^N_0 \in \PPP_{\!sym}(E^N)$, with the additional assumption supp$\, G^N_0 \subset \BB\SS^N$ in the {\bf HS} case, and $G^N(t)$ the solution to 
the Boltzmann-Kac system \eqref{eq:BKs}. 
Consider  $f_0 \in \PPP(E)$ and $f(t)$ the corresponding solution to the Boltzmann equation  \eqref{eq:Beq}. 

(a) If $G^N_0$ is $f_0$-chaotic, then $G^N(t)$ is $f(t)$-chaotic.

(b) Better, for the  {\bf MG} case and if  $G^N_0 = f_0^{\otimes N}$,  then 
$$
  \qquad \sup_{t \in [0,T]} W_1(G^N_j(t),f(t)^{\otimes j})  \le {C_{j,T} \over N} \quad
  \hbox{for any fixed} \, j \ge 1.
 $$
\end{theo}

In the above statement, we  define the $j$-th marginal $G^N_j \in \PPP(E^j)$, $1 \le j \le N$,  of $G^N$  by  
\beqn\label{def:GNj}
G^N_j = \int_{E^{N-j}} G^N dz_{j+1} ... dz_N. 
\eeqn
Moreover,  for $F,G \in \PPP(E^j)$ we define by $W_1(F,G)$ the (renormalized) Monge-Kantorovich-Wassertstein (MKW) 
distance
\beqn\label{def:W1PPPE}
W_1(F,G) := \inf_{\pi \in \Pi(F,G)} \int_{E^j \times E^j} \Bigl( {1 \over j} \sum_{i=1}^j |x_i-x_j| \wedge 1 \Bigr) \, \pi(dX,dY),
\eeqn
where $\Pi(F,G)$ stands for the set of probabilty measures on $E^j \times E^j$ with given first marginal $F$ and second  marginal $G$. 

\smallskip
It is worth emphasizing again that proving the propagation of chaos result as formulated in Theorem~\ref{theo:Kac1}
imply that we are able to identify  the large number
of particles limit (law of large numbers in the mean field limit) of the system of particles, or in other words, we have derived the (space homogeneous)
Boltzmann equation from a microscopic description (of the physical system). 

\smallskip
For the {\bf MG} model, 
the propagation of chaos (without rate and next with rate) has been proved by Kac \cite{Kac1956}, McKean \cite{McKean1967,McKean-TCL}, Grunbaum \cite{Grunbaum}, Tanaka \cite{T2}, Graham, M\'el\'eard \cite{GM}
using (except in \cite{Grunbaum}) some tree arguments (Wild sum, stochastic tree). These kind of arguments are very specific to the {\bf MG}
 model. 


For the {\bf  HS} model, the propagation of chaos result (without rate) has been proved by Sznitman (1984, \cite{S1}) using a nonlinear Martingale approach, some compactness of the system and uniqueness of the limit arguments. An alternative proof is suggested in  Arkeryd et al (1991, \cite{ArkerydCI99}) 
following a ``BBGKY hierarchy" approach. 

 \smallskip
For the {\bf  HS} model again, in order to be able to apply the first part of Theorem~\ref{theo:Kac1}, we have to build a sequence of initial data $G_0^N$ which satisfies both
properties
$$
\hbox{supp}\, F^N_0 \subset \KK\SS^N \,(\hbox{or}\, \BB\SS^N )
\quad\hbox{and}\quad
F^N_0 \,\, \hbox{is}\,\, f_0\hbox{-chaotic.}
$$
A first answer to that issue is the following.  

\begin{theo}[Kac \cite{Kac1956}; Carlen et al. \cite{CCLLV}] \label{theo:Kac2}
Consider $f_0 \in L^1_4 (E) \cap L^p(E)$, $p > 1$, $E=\R$.  There exists $F^N_0 \in \PPP(E^N)$ such that 

(a) $\hbox{supp}\, F^N_0 \subset \KK\SS^N$;

(b)  $F^N_0$ is $f_0$-chaotic;

(c) $F^N_0$ is $f_0$-entropy chaotic. 

\end{theo}

\smallskip\noindent
Here, we say that a sequence $(G^N)$ of $\Ps(E^N)$ is $f$-entropy chaotic if  
\beqn\label{def:entropichaos}
(G^N) \hbox{ is } f\hbox{-(Kac's) chaotic and } H(G^N|\gamma^N) \to H(f|\gamma),
\eeqn
where 
$$
H(G^N|\gamma^N) := \frac1N \int_{\KK\SS^N} G^N \log { dG^N \over d\gamma^N}, 
\quad
H(g|\gamma) := \int_{E} g \, \log {g \over \gamma}, 
$$
and $\gamma$ is the normalized gaussian function.

\subsection{Problem 2: Known results on the convergence to the equilibrium uniformly with respect to the number of particles}

Kac believed that one can obtain a rate of convergence in the large time asymptotic for the {\bf nonlinear} Boltzmann equation
from the same result for the {\bf linear} Boltzmann-Kac system in large increasing dimension. 
This has motivated beautiful works on the ``Kac spectral gap problem'', i.e. the study of this relaxation rate in a $L^2$ setting,
for the Kac's $N$-particle system first  \cite{Kac1956,janvresse,maslen,CCL2003,CCL-preprint} and next for
the Boltzmann-Kac system  \cite{CCL2003,CarlenGeronimoLoss2008}. 

\begin{theo}\label{theo:SpectralGap}
For both {\bf M} and 
{\bf HS} models, there exists $ \delta > 0$ such that for any $N \ge 1$
$$
\Delta_N := \inf  \{ - \langle h, \Lambda^N h \rangle_{L^2}, \,\, \langle h, 1 \rangle_{L^2} = 0, \,\, \| h \|_{L^2}^2=1 \} \ge \delta  ,
$$
where $\langle \cdot, \cdot \rangle_{L^2}$ and $\| \cdot \|_{L^2}$ stand for the scalar product and the norm of  $L^2(\BB\SS^N; d\gamma^N)$. 

That spectral gap estimate implies that for any $G^N_0 = h_0 \, \gamma^N \in \PPP_{\! sym} (E^N)$, $h_0 \in L^2$, the associated solution $G^N$ to the Boltzmann-Kac system  \eqref{eq:BKs} can be written as $G^N = h(t) \, \gamma^N$ and satisfies 
\beqn\label{eq:hNto1}
 \| h^N(t) - 1 \|_{L^2} \le e^{-\delta \, t} \, \| h^N_0 - 1 \|_{L^2}.
\eeqn
\end{theo}

\smallskip

Some few remarks are in order.

\smallskip
(a) Theorem~\ref{theo:SpectralGap} does not answer Problem 2 because if $G^N_0= h^N_0 \, \gamma^N$ is $f_0$-chaotic then $\| h^N_0 - 1  \|_{L^2} \ge A^N$, with $A > 1$,  and we have to wait a time proportional to  $N$ in order that \eqref{eq:hNto1} implies any convergence to the equilibrium.

\smallskip
(b) The entropy fits better for a $N\to\infty$ asymptotic. However in that case, the ``spectral" gap  
$$
\Delta'_N := \inf  \{ - \langle (\log dG/d\gamma^N)/N, \Lambda^N G \rangle / H(G|\gamma^N) \} ,
$$
satisfies  $\Delta'_N \ge 1/N$ (Villani \cite{VillaniCerConjecture}) and  $\limsup \Delta'_N = 0$ (Carlen et al \cite{CCL2003}), and that cannot answer Problem 2 either. 
 
\smallskip
(c) On the other hand, the exponential rate of convergence to the equilibrium of the solutions to the nonlinear Boltzmann equation 
has been established by a direct PDE approach :  
\beqn\label{eq:ftogamma}
 f_0 \in \PPP_2(E) \qquad D(f(t), \gamma) \le C_{f_0} \, e^{-\lambda \, t}
\eeqn
for an appropriate distance $D$ on $\PPP(E)$ and where $\gamma$  is the Gaussian function associated to the initial datum $f_0$. 
There is a so huge number of works on that topics that we cannot quote all of them. Let us just say that the story began with the work by T. Carleman 
(1933, \cite{MR1555365}) and we refer to \cite{LuMou} and the references therein for the {\bf HS}  model and to  \cite{CGT}
and the references therein for the {\bf M}  model. 


\bigskip
\section{A reverse answer to Kac's program }
\label{sec:KacProgram}
\setcounter{equation}{0}
\setcounter{theo}{0}

\subsection{Our contributions to Kac's program  } We give some possible answers to the three problems formulated by Kac.

\medskip
$\bullet$ In collaboration with C. Mouhot and B. Wennberg in \cite{MMWchaos}  we develop a new quantified method
for proving propagation of chaos inspired from Grunbaum's work \cite{Grunbaum}. We illustrate our approach on several
$N$-particle system including the Vlasov, the McKean-Vlasov and the Boltzmann models under mild assumptions
on the coefficients.  


\medskip
$\bullet$  In collaboration with M. Hauray in \cite{HaurayMischler}, and in a next  work \cite{KleberSphere} by K. Carrapatoso, we revisit the notion of chaos by Kac on a flat space $E^N$ and we establish some links with the notion of entropic chaos and Fisher chaos (that we introduce). 
We then give several extensions to the Kac's spheres $\KK\SS^N$ framework, to the Boltzmann's spheres $\BB\SS^N$ framework and to the
De Finetti, Hewitt \& Savage's mixture (without chaos) framework. 
 
\medskip
$\bullet$ In collaboration with C. Mouhot in \cite{MMinvent}, we establish a quantitative propagation of chaos estimate  for the  realistic Boltzamnn-Kac system associated to hard spheres interactions and to true Maxwell molecules interactions, and then
we make significant progress in the answer to  Kac's problem 1, since we improve  Sznitman's result  \cite{S1} on the propagation of chaos for the hard spheres models (without quantitative estimate) as well as the many works (Kac \cite{Kac1956},  McKean \cite{McKean1967}, Grunbaum \cite{Grunbaum}, Tanaka \cite{T2}, Graham and M\'el\'eard~\cite{GM}, Fournier and M\'el\'eard~\cite{FM7,FM10}, Kolokoltsov \cite{Kolokoltsov}, Peyre \cite{Peyre})  which deal with the Maxwell molecules  model with Grad cutoff. 

\smallskip
$\rhd$ Our estimate is uniform in time  and  that make possible to answer Problem 2 by Kac on the convergence to equilibrium with uniform rate  with respect to the number of particles but in the reverse sense that which Kac imagined. Indeed, the later long time convergence is deduced from the uniform in time propagation of chaos estimate together with the known results about the (not uniform in $N$) convergence to the uniform density $\gamma_N$ (Theorem~\ref{theo:SpectralGap}) and the known result about convergence to the (Maxwellian) equilibrium for the nonlinear Boltzmann equation (estimate \eqref{eq:ftogamma}). 

 \smallskip
$\rhd$ We prove the entropic chaos of the Boltzmann-Kac system and that is the first derivation of the $H$-Theorem for the Boltzmann equation
starting  from a microscopic description of the system, providing a possible answer to Problem 3 by Kac.

\medskip
Our approach developed in  \cite{MMWchaos,MMinvent} gives an alternative method to the classical coupling method initiate by A. Sznitman \cite{S6}.
For that later, we refer to the work~\cite{MR2731396} and the reference therein for recent development of coupling method
for the McKean-Vlasov model, as well as  to the work in collaboration with N. Fournier \cite{FM} and the references therein 
for the application of the coupling method to a Kac-Boltzmann related (but simpler) collisional system.
 It is worth mentioning that the estimates obtained in  \cite{FM}  are more accurate (in the number of particles) than those 
 obtained thanks to the method developed in  \cite{MMWchaos,MMinvent}, but they are    however local in time. 

\medskip
In the next section we give a more detailed statement of our results. 

\subsection{Propagation of chaos for the Boltzmann-Kac  system  }

\begin{theo}[{\bf M} case, \cite{MMinvent}] \label{Th:MM-MM} For any $f_0 \in \PPP(E)$,  which satisfies additional smoothness conditions,  there exists a sequence $G^N_0$ 
of  $f_0$-chaotic initial data such that 
$$
\forall \, N \ge 1 \qquad  \sup_{t \ge 0} W_1(G^N(t),f(t)^{\otimes N}) \le {C \,   \over N^\bullet},  
 $$ 
 for some positive (small!) exponent $\bullet > 0$, and 
\beqn\label{eq:boundFisher}
\forall \, N \ge 1 \qquad  \sup_{t \ge 0} I(G^N(t)|\gamma^N) \le C,  
\eeqn
 where the functional $I$ stands for the relative Fisher information (see below). 
\end{theo}

\medskip
Our result is mainly based on : 

\smallskip
$-$ an appropriate  differential calculus in the probability measures space $\PPP(E)$;

\smallskip
$-$ some accurate stability estimates on the nonlinear Boltzmann equation in $\PPP(E)$;

\smallskip
$-$ the equivalence of different ways to measure the chaos;

\smallskip
$-$ a functional law of large number estimate on the initial chaos. 

\medskip
We deduce then  from Theorem~\ref{Th:MM-MM} that   
$$
\forall \, N \ge 1 \qquad  \sup_{t \ge 0} \Bigl(  |H(G^N(t)|\gamma^N) - H(f(t)|\gamma)| +  \| G^N_j (t) - f^{\otimes j} \|_{L^1} \Bigr)\le   {C_j \over N^\bullet}, 
 $$
\beqn\label{eq:MM-GNtogammaN}
\forall \, t > 0 \qquad  \sup_{N \ge 1} \Bigl(  W_1(G^N(t),\gamma^N) +  H(G^N(t)|\gamma^N) +  \| G^N_j (t) - \gamma^{\otimes j} \|_{L^1}\Bigr) \le  {C_j \over t^\bullet} . 
\eeqn
Again, these estimates give an answer to the three problems of  Kac's program. 

\bigskip
For the Hard spheres model we have the following (weaker) variant: 
 
\begin{theo} [{\bf HS} case, \cite{MMinvent}] \label{Th:MM-HS}
For any $f_0 \in \PPP(E)$,  which satisfies additional smoothness conditions, there exists a sequence $G^N_0$ 
of  $f_0$-chaotic initial data such that 
$$
\forall \, N \ge 2 \qquad  \sup_{t \ge 0} W_1(G^N(t),f(t)^{\otimes N}) \le {C \,   \over (\log N)^\bullet }  
 $$
for some positive  exponent $\bullet > 0$,  and $H(G^N(t)|\gamma^N)  \to H(f(t)|\gamma)$ for any $t \ge 0$ (non uniformly in time).
 \end{theo}
 
We also deduce 
\beqn\label{eq:HS-GNtogammaN}
\forall \, t \ge 2 \qquad \sup_{N \ge 1} W_1 (G^N(t),\gamma^N) \le {C  \over (\log t)^\bullet}. 
\eeqn

\bigskip
\section{Uniformly in time chaos estimate}
\label{sec:KacProgram}
\setcounter{equation}{0}
\setcounter{theo}{0}

\subsection{ Weak uniform in time quantitative chaos propagation}

We give the cornerstone estimate of the quantitative propagation of chaos method developed in \cite{MMWchaos,MMinvent} for which we 
present next a sketch of the proof. 

\begin{theo} [{\bf HS} \& {\bf M}, \cite{MMinvent}]  \label{theo:Kac3} Under the assumptions and notations of Theorems~\ref{Th:MM-MM} \& \ref{Th:MM-HS}, there holds for any $\eps > 0$ and for some suitable modulus of continuity $\Theta$ 
 \[
 \forall \, N \ge 1 \qquad \sup_{t \ge 0} \| G^N_2 - f^{\otimes 2} \|_{{\FF}'} \le {C_\eps  \over N^{1-\eps}} +\Theta ( \langle G^N_0, W_1(\mu^N_V,f_0) \rangle ) 
\]   
in a {weak dual norm} $ \| \cdot  \|_{{\FF}'}$, with  ${\FF}$ a space of smooth functions $ \subset UC_b(E^2)$.
\end{theo}

The difficulty here is to compare the two solutions  $G^N_t$ on $\PPP(E^N)$ and $f_t$ on $\PPP(E)$ which do not belong to the same functional space.  The idea is to compare the {dynamics } associated to \eqref{eq:BKs} and to \eqref{eq:Beq}  {(and not only the two solutions $G^N_t$ and $f_t$)}
{ both in the same space $C_b(P(E))$} through some relevant  {``projections"} : 
 
\begin{itemize}

\item we project the $N$-particle dynamics thanks to the empirical measures map $E^N \to P(E)$, 
$V \mapsto \mu^N_V$, with 
$\displaystyle{ \mu^N_V(dz) := {1 \over N} \sum_{i=1}^N \delta_{v_i}(dz)}$. 

\item we project the mean-field dynamics by pullback.

\end{itemize}

Similar ideas have been used by  Grunbaum in \cite{Grunbaum}
(which is a main source of inspiration) and by Kolokoltsov in \cite{Kolokoltsov} (independently).

\subsection{ Sketch of the proof : splitting and estimate the three terms separately   }
\vspace{-0.0cm} 

We split 

\bean
 \left \langle G^N_t -  {f^{\otimes N}_t} , \varphi \otimes 1^{\otimes N-k} \right\rangle
 &=&  
 \left \langle G^N_t , \varphi \otimes 1^{\otimes N-k} -  {R_\varphi }( \mu^N_V)  \right\rangle  \qquad (= T_1)
 \\
&&+  \left \langle G^N_t,  {R_\varphi }( \mu^N_V)  \right\rangle -  \left \langle G^N_0,   {R_\varphi } ( {S^{NL}_t }\mu^N_V) \right\rangle
 \qquad (= T_2)
\\
&&+  \left \langle G^N_0,   {R_\varphi } ( {S^{NL}_t } \mu^N_V) \right\rangle - \left \langle  {f_t^{\otimes k}}, \varphi \right\rangle
 \qquad (= T_3) , 
\eean
where $  {R_\varphi }$ is the ``polynomial function" on $\PPP(\R^3)$ defined by 
$$
  {R_\varphi }(\rho) = \int_{E^k} \varphi \, \rho(dv_1) \, ... \, \rho(dv_k)
$$
and we recall that $ {S^{NL}_t }$ is the nonlinear semigroup defined by \eqref{eq:defSNLt}.

\subsection{ The term $T_2$. } We rewrite the term $T_2$ as    
\bean
T_2 
&:=&   \left \langle G^N_t, R_\varphi ( \mu^N_V)  \right\rangle -  \left \langle G^N_0, R_\varphi (S^{NL}_t \mu^N_V) \right\rangle 
\\
&=& 
 \left \langle G^N_0,  T^N_t (R_\varphi \circ \mu^N_V)  -  (T_t^\infty  R_\varphi ) ( \mu^N_V) \right\rangle, 
\eean
with
\begin{itemize}
\item  {$T^N_t = $} the dual semigroup (acting on  {$ C_b(E^N)$}) of the N-particle flow $G^N_0 \mapsto G^N_t$, 
\item  {$T^\infty_t = $} the pushforward semigroup  (acting on  {$C_b(\PPP(E))$}) of the nonlinear semigroup $S^{NL}_t$
defined by $(T^\infty \Phi)(\rho) := \Phi(S^{NL}_t \rho)$. 
 \end{itemize} 

\medskip\noindent
Introducing 
\begin{itemize}
\item  {$\pi_N = $} the projection from   {$C_b(\PPP(E))$} onto   {$C_b(E^N)$} defined by $(\pi_N \Phi)(V) = \Phi(\mu^N_V)$ 
for any $\Phi \in C_b(\PPP(E))$ and $V \in E^N$, 
\end{itemize} 
we may rewrite $T_2$ as  the difference of the two dynamics in $C_b(P(E))$  in the following way
\bean
T_2 
 &=& 
 \left \langle G^N_0,  (T^N_t \pi_N - \pi_N T_t^\infty ) \, R_\varphi  \right\rangle.
\eean

\medskip\noindent
 Thanks then to Trotter-Kato formula, we have 
\bean
T_2  
&=&  \left \langle G^N_0, \int_0^T  T^N_{t-s}  \, (\Lambda^N \pi_N - \pi_N \Lambda^\infty) \, T_s^\infty   \, ds \, R_\varphi  \right\rangle
 \\
&=&  \int_0^T  \Bigl \langle  {\underbrace{ G^N_{t-s}}_{moment}},   { \underbrace{(\Lambda^N \pi_N - \pi_N  \Lambda^\infty )}_{consistency} } \, { \underbrace{T_s^\infty \Phi}_{stability} }    \Bigr\rangle \, ds ,  
\eean
with  $\Phi := R_\varphi$.

 \bigskip
In order to get a bound on $T_2$,  we have to verify the four following assertions: 
\begin{itemize}
\item  {(A1)} $G^N_t$ has enough  {polynomial} bounded moments. 
\item   {(A2)} $\Lambda^\infty \Phi (\rho) = \langle Q(\rho),D\Phi (\rho) \rangle$.
\item   Consistency result {(A3)}: the difference of generator applied on "smooth" functions is of order $1/N$ in the sense 
$$
(\Lambda^N \pi_N \Phi) (V) = \langle Q(\mu^N_V),D\Phi (\mu^N_V) \rangle + \OO ( [\Phi]_{C^{1,\eta}} / N).
$$ 
Here a ``smooth" function $\Phi$ means that we can perform an  expansion of $\Phi$ up to order $1+\eta$ in each point  of  $\PPP(E)$ seen as an embedded manifold  
of ${\FF}'$  (that is a much simpler notion than  the ``differential calculus" developed in the ``gradient flow theory"). 

\item  Last, we have to check that $\Phi_s := T^\infty_s \Phi $ remains a ``smooth" function, and that is a consequence of 
the stability result  {(A4)}:  
$S^{NL}_t \in C^{1,\eta} (\PPP(E); \PPP(E))$. 
\end{itemize} 

Among these assertions (A1)-(A4), the probably newer and more technical one is the last one that we state now (in the {\bf HS} case).

\begin{theo} [\cite{MMinvent}] \label{theo:SNLholder}
For any $f_0,g_0 \in P(E)$ with energy $1$ and mean velocity $0$, there holds
\[
\forall \, t \ge 0 \qquad  \| f_t - g_t - h_t \|_{L^1}  \le C \, M_k(f_0+g_0) \, e^{-\lambda \, t} \, \| f_0 - g_0 \|_{L^1}^{1+\eta},
\]    
where $h_t := $ solution to the linearized Boltzmann equation around $f_t$ and $k > 2$. 
\end{theo}

\noindent
Theorem~\ref{theo:SNLholder} is nothing but a refinement version of uniqueness and moments estimates obtained by many authors, among who one can quote:  Povzner, Arkeryd, Elmroth, Desvillettes, Carlen, Carvalho, Toscani, Gabetta, Villani, Lu, Wennberg, M., Mouhot. 

 \smallskip
 In a simple case (which corresponds to what one can do locally in time for the {\bf MG} model), we then may bound $T_2$ 
 in the following way:
\bean
T_2 &\le&  \int_0^T  M_0(G^N_{t-s}) \,  \| (\Lambda^N \pi_N - \pi_N \Lambda^\infty) \, (T_s^\infty  \, R_\varphi) \|_{L^\infty(E^N)}  \, ds 
\\
&\le& \int_0^T {C \over N}  \, [ T_s^\infty  \, R_\varphi ]_{C^{1,\eta}}  \, ds  \qquad\quad (\hbox{using } (A1), (A2), (A3))
\\
&\le& {C\over N}  \int_0^T   \, [ R_\varphi \circ S^{NL}_t  ]_{C^{1,\eta}}  \, ds 
\\
&\le& {C\over N}  \int_0^T [ R_\varphi ]_{C^{1,1}} \, [ S^{NL}_t  ]_{C^{1,\eta}}  \, ds 
\\
&\le& {C \over N} \, k^2  \, \| \varphi \|_{W^{2,\infty}}  \int_0^T  [S^{NL}_t  ]_{C^{1,\eta}}  \, ds 
\\
&\le&  {C \over N} \, k^2  \, \| \varphi \|_{W^{2,\infty}} 
\qquad\quad\qquad\quad\qquad\quad (\hbox{using } (A4)). 
\eean

\subsection{ Estimate of $T_1$. } Thanks to a combinatory trick which is due to Grunbaum \cite{Grunbaum}, 
we find
\bean
|T_1| 
&=&  \left |  \left \langle G^N_t , \varphi \otimes 1^{\otimes (N-k)} (V) - R_\varphi ( \mu^N_V)   \right\rangle \right|
\\
&=&  \left |  \left \langle G^N_t ,   { \widetilde{\varphi \otimes 1^{\otimes (N-k)}} (V) - R_\varphi ( \mu^N_V) }  \right\rangle \right|
\\
&\le&
   \left \langle G^N_t ,    { {2 \, k^2  \over N} \, \| \varphi \|_{L^\infty(E^k)}  }   \right\rangle 
\\
&=&  {2 \, k^2  \over N} \, \| \varphi \|_{L^\infty(E^k)} , 
\eean
where we use that  $G^N$ is symmetric and a probability and we introduce the symmetrization function associated to $\varphi \otimes 1^{\otimes (N-k)}$ by 
$$
 \widetilde{\varphi \otimes 1^{\otimes (N-k)}} (V) = {1 \over \sharp \SIGMA_N} \sum_{\sigma \in \SIGMA_N} \varphi \otimes 1^{\otimes (N-k)}(V_\sigma).
$$

\subsection{ Estimate of $T_3$. } We claim  that the nonlinear flow $S^{NL}_t$  is $\Theta$-Holder continuous in the sense that 
$$
 (A5) \qquad W_1(f_t,g_t) \le  {\Theta(} W_1(f_0,g_0) { )}\quad \forall \, f_0,g_0 \in \PPP_{ {exp}}(E)
$$
for $\Theta(u) = u$ ({\bf M} model) and $\Theta(u) = |\log|u|\wedge1|^{-1}$ ({\bf HS} model). Such an estimate has been proved by 
Tanaka \cite{T1} and Toscani-Villani \cite{TV} in the case {\bf M} and it is mainly a  consequence of Fournier-Mouhot \cite{Fo-Mo}
in the case of {\bf HS}. As a consequence, and similarly as  for the empirical measures method, we find
\bean
|T_3| &=& \left|  \left \langle G^N_0, R_\varphi (S^{NL}_t \mu^N_V) \right\rangle - \left \langle (S^{NL}_t f_0)^{\otimes k}, \varphi \right\rangle \right|
\\
&=& \left|  \left \langle G^N_0, R_\varphi (S^{NL}_t \mu^N_V) - R_\varphi (S^{NL}_t f_0) \right\rangle \right|
\\
&\le&  { [R_\varphi]_{C^{0,1}} } \left \langle G^N_0,  { W_1 (S^{NL}_t \mu^N_V , S^{NL}_t f_0)}  \right\rangle 
\\
&\le& {  k \, \| \nabla \varphi \|_{L^\infty(E^k)} } \,   \left \langle G^N_0,   { \Theta (W_1 (\mu^N_V ,  f_0) ) } \right\rangle
\\
&\le&  k \, \| \nabla \varphi \|_{L^\infty(E)} \,   {\Theta \left(   \Black{   \left \langle G^N_0,   W_1 (\mu^N_V ,  f_0)    \right\rangle } \right)}
\eean
where 
$$
 {  [R_\varphi]_{C^{0,1}} := \sup_{W_1(\rho,\eta) \le 1} |R_\varphi(\eta) - R_\varphi(\rho)| = k \, \| \nabla \varphi \|_{L^\infty} } .
$$
 
\bigskip

\section{Kac's chaos and related problems}
 
 \label{sec:LongTimeBehave}
\setcounter{equation}{0}
\setcounter{theo}{0}

We present in this section several results obtained in \cite{MMinvent,HaurayMischler,KleberSphere} 
about Kac's chaos. They will be useful in order to deduce Theorems~\ref{Th:MM-MM} and \ref{Th:MM-HS} from 
Theorem~\ref{theo:Kac3}. 

\smallskip
More precisely, for the flat spaces $E^N$, we establish

- the equivalence between the several ways to measure Kac's chaos;

- some links between Kac's chaos and stronger definitions of chaos (in the sense of entropy and Fisher information); 

- some rate of convergence in the functional law of large numbers. 

\smallskip
Next, we generalize all these results to the Kac's spheres and Boltzmann's spheres framework. 

\smallskip
It is worth emphasizing that these results can also been adapted to a situation without chaos, that is for a sequence of
$\PPP(E^N)$ which converges to some mixture measure in the sense of De Finetti, Hewitt and Savage.
However, we will not need that extensions in the present context, and we refer the interested reader to \cite{HaurayMischler}
as well as to \cite{FHM} for an application of that tools.

\subsection{Several definitions of Kac's chaos }

Kac's chaos has be formulated by Kac, Grunbaum and Sznitman in several ways that we recall here. 
Again, Kac's chaos aims to formalize the intuitive idea for
a family of exchangeable stochastic variables with values in $E^N$ to be asymptotically independent. Using the diagram
\bean
\ZZ^N = (\ZZ^N_1, ..., \ZZ^N_N) \in E^N &\rightarrow& G^N := \LL(\ZZ^N) \in \PPP_{\! sym}(E^N)
\\
 \updownarrow \qquad \qquad  \quad && \qquad  \qquad \updownarrow
\\
\mu^N_{\ZZ^N} := {1 \over N} \sum_{i=1}^N \delta_{\ZZ^N_i} \in \PPP(E)  
&\rightarrow&  \ \hat G^N := \LL(\mu^N_{\ZZ^N}) \in \PPP(\PPP(E)),
\eean
Kac's chaos can be formulated in the two following equivalent ways 
 
$-$ $G^N_j  \wto f^{\otimes j}$ as $N\to\infty$, for some (any) fixed $j \ge 2$ ($j \ge 1$),

$-$ or $\mu^N_{\ZZ^N} \Rightarrow f$  as $N\to\infty$ (in law for stochastic variables with values in $\PPP(E)$),

$-$ and it is known to be implied by $G^N \sim f^{\otimes N}$ asymptotically as  $N\to\infty$.

\medskip\noindent
In a more precise (but equivalent) way $G^N$ is  $f$-(Kac's) chaotic if :

\medskip
$-$ $W_1(G^N_j, f^{\otimes j} ) \to 0$ as $N\to\infty$, for some (any) fixed  $j \ge 2$ ($j \ge 1$), 

\smallskip
$-$ or $\EEE ( W_1(\mu^N_{\ZZ^N}, f)) \to 0$ as  $N\to\infty$,

\smallskip
$-$ and both are implied by $W_1(G^N,f^{\otimes N})  \to 0$ as $N\to\infty$.

\subsection{Equivalence between the measures of Kac's chaos }

For $G^N \in \PPP(E^N)$, we define  $\hat G^N \in \PPP(\PPP(E))$ by 
$$
\forall \, \Phi \in C_b(\PPP(E))  \qquad \langle \hat G^N,\Phi \rangle := \langle G^N,\pi_N \Phi\rangle
$$
and then  $\WW_{W_1}$ as the MKW distance in $\PPP(\PPP(E))$ associated to the $W_1$ distance in  $\PPP(E)$.

\begin{theo}[\cite{HaurayMischler}]  \label{theo:EquivKacChaos}
In $E:=\R^d$, the following convergences are equivalent : 
\bear 
\label{eq:Omegaj}
&& \Omega_{j}(G^N;f):=W_1(G^N_j, f^{\otimes j} ) \to 0 \ \hbox{ as } \ N\to\infty, \hbox{ for any fixed  } j \ge 2,
\\ 
\label{eq:Omegainfty}
&& \Omega_{\infty}(G^N;f):=W_{W_1}(\hat G^N, \delta_f)  = \EEE ( W_1(\mu^N_{\ZZ}, f)) \to 0  \ \hbox{ as } \  N\to\infty,
\\
\label{eq:OmegaN}
&& \Omega_N(G^N;f):= W_1(G^N, f^{\otimes N} ) \to 0  \ \hbox{ as } \  N\to\infty.
\eear
More precisely, for any $M, \, k > 1$, there exist  $\alpha, \, C  > 0 $ such that for any 
$ f \in \PPP (E)$, for any $G^N \in \PPP_{sym} (E^N)$ such that  $M_k(G^N_1), M_k(f) \le M$, there holds 
\[
\forall \,  j , \ell \in \{2, ..., N,\infty \}  
\qquad 
\Omega_j \le C \, \left( \Omega_\ell^{\alpha_1} + {1 \over N^{\alpha_2}} \right). \qquad
\]     
\end{theo}

Let us emphasize that Theorem~\ref{theo:EquivKacChaos} gives the equivalence of \eqref{eq:Omegaj} and  \eqref{eq:Omegainfty}
in a quantitative way, and it also establishes the new assertion that \eqref{eq:Omegaj} implies \eqref{eq:OmegaN}. 

\medskip
We just give some ideas about the proof. 

\smallskip
 $\rhd$  On the one hand, taking $F^N = f^{\otimes N},G^N \in \PPP_{\! sym}(E^N)$, or equivalently in probabilistic words, taking $\ZZ$ of law $G^N$ and $\YY$ of law $f^{\otimes N}$, we have 
$$
W_1(G^N,F^N) = \EEE ( W_1 (\mu^N_\ZZ, \mu^N_\YY)) =  \WW_{W_1} (\hat G^N,\hat F^N) \simeq \WW_{W_1} (\hat G^N, \delta_f) +  {  \WW_{W_1 } (\hat F^N,\delta_f) },
$$
where the first identity is proven thanks to an optimal mass transportation argument and the last equivalence is nothing but the triangular inequality. 
That establishes the equivalence between $\Omega_N$ and $\Omega_\infty$.

\smallskip

 $\rhd$ On the other hand, we just have to estimate $\Omega_\infty$ by $\Omega_2$ since the bound of $\Omega_2$ by $\Omega_N$ is trivial. 
 Using an interpolation argument, we may replace the $W_1$ distance in $\PPP(E)$ by 
 a Hilbert norm  $\| \cdot \|_{H^{-s}}$, $s > d/2$, which enjoys better algebraic properties, and next just develop as usually, in the following way 
\bean
C^{-1} \,  \WW_{W_1} (\hat G^N, \delta_f)^\bullet &=& C^{-1} \, \EEE ( W_1 (\mu^N_\ZZ, f))^\bullet 
\\
&\le& \EEE (  \| \mu^N_\ZZ -f \|^2_{H^{-s}} )  \le \, C_s \, W_1(G^N_2,f^{\otimes 2})+ {C_s \over N}
 \eean
 for some exponent $\bullet > 1$.

\subsection{ Chaos by Kac, Boltzmann and  Fisher}
\medskip
We carry on the work \cite{CCLLV} by Carlen, Carvalho, Le Roux, Loss and Villani on the entropic (or Boltzmann's) chaos. 
Following \cite{CCLLV} we may define the entropic chaos by \eqref{def:entropichaos}. 

\smallskip
Next, considering a sequence $(G^N)$ of  $\PPP_{\!sym}(E^N)$ and $f \in \PPP(E)$, we say that  $(G^N)$ is $f$-Fisher's chaotic if  $G^N_1 \ \wto \ f$ and $I(G^N) \to I(f)  < \infty$, where 
$$
 I(F) := {1 \over j} \int_{E^j} {|\nabla F|^2 \over F} 
\quad\hbox{for}\quad
F \in \PPP(E^j).
$$
 
 \begin{theo}  [\cite{HaurayMischler}] \label{theo:FBKchaos} Consider a sequence $(G^N)$ of $\PPP_{\!sym}(E^N)$ such that $M_k(G^N_1)$ is bounded, $k > 2$, 
 and $f \in \PPP(E)$. 
 In the list of assertions below, each one implies the assertion which follows

\medskip

$(a)$ $(G^N)$ is  $f$-Fisher's chaotic;

$(b)$  {$(G^N)$ is   $f$-Kac's chaotic and  $I(G^N)$ is bounded;}

$(c)$  {$(G^N)$ is  $f$-Boltzmann's chaotic;}

$(d)$ $(G^N)$ is  $f$-Kac's chaotic.

\end{theo}

The most unexpected  and interesting implication is maybe (b) $\Rightarrow$ (c). It is a simple consequence of the HWI inequality 
$$
H(f_2|\rho) \le H(f_1|\rho) + \sqrt{I(f_0|\rho)} \, W_2(f_0,f_1)
$$
of  Otto and Villani \cite{OttoVillani},  which is itself a kind of (dimensionless) variant of the wellkonwn (and mere consequence of 
the Cauchy-Schwarz inequality) interpolation inequality
$$
\| g \|_{L^2} \le \| g \|_{H^1}^{1/2} \, \| g \|_{H^{-1}}^{1/2}.
$$
Indeed, using twice the HWI inequality, we have  
$$
\Bigl|H(G^N | \gamma^{\otimes N}) - H(F^N | \gamma^{\otimes N}) \Bigr| \le  \Bigl( \sqrt{ I(F^N | \gamma^{\otimes N}) }  +  \sqrt{ I(G^N | \gamma^{\otimes N}) } \Bigr) \, W_2 (G^N,F^N)
$$
and  (b) $\Rightarrow$ (c) follows thanks to Theorem~\ref{theo:EquivKacChaos}. 
The implications  (a) $\Rightarrow$ (b) and  (c) $\Rightarrow$ (d) can be deduced from the lsc property
$$
H(G_2) \le \liminf H(G^N) \quad\hbox{if} \quad G^N_2 \wto G_2
$$
and from the fact that the saturation equalities $H(G_2) = H(f \otimes f)$,  $H(G_1) = H(f)$ (of a related convexity inequality) implies $G_2 = f \otimes f$.

\subsection{ Functional law of large numbers }
 
There exist many works, for instance by Glivenko, Cantelli, Talagrand, Rachev-R\"uschendorf, Boissard-LeGouic, Barthe-Bordenave,
on the functional law of large numbers, that is on the convergence (or for the more recent works on the proof of a rate for the convergence) of the empirical measures 
$$
\mu^N_{\ZZ^{N}} \Rightarrow f \quad\hbox{as}\quad N\to\infty, 
$$ 
when  $\ZZ^{N} \sim f^{\otimes N}$.  

We give an accurate (and sometimes almost optimal) version of the above convergence result that we also generalize to some situations  when
we only assume that $\ZZ^{N}$ is $f$-Kac's chaotic. 

\begin{theo} [\cite{MMinvent,HaurayMischler,KleberSphere}] 
Under suitable boundedness assumptions, there holds
\beqn\label{eq:functLLN}
 \Omega_{\infty}(G^N;f):=W_{W_1}(\hat G^N, \delta_f) = \EEE(W_{1}(\mu^N_{\ZZ^N} ,f) )  \le C/N^{\alpha}
\eeqn
for 
\begin{itemize}
\item[(i)]  any  $\alpha < \alpha_{c} \approx 1/d$ when $G^N = f^{\otimes N}$;

\item[(ii)] any  $\alpha < 1/2$ when  $G^{N } = \gamma^N$ and $f = \gamma =$ the Gaussian function, as well as 
\beqn\label{eq:PoincarLemma}
 \Omega_{N} (\gamma^N;\gamma) := W_1(\gamma^N,\gamma^{\otimes N} ) \le C/N^{1/2};
\eeqn
\item[(iii)] any  $\alpha < 1/(2+2/k)$ when  $G^{N } = $ the ``conditioned to the Kac's spheres of the product measures $f^{\otimes N}$", $f \in \PPP_{k}(\R)$, 
and a similar result replacing the ``Kac's spheres" by the ``Boltzmann's spheres". 

\end{itemize}    

 \end{theo}

\noindent
Let us make some comments:

$\rhd$ In the case (i) and $f$ has compact support, the estimate \eqref{eq:functLLN} with  $\alpha = \alpha_{c} $  is true, see \cite{BoissardLeGouic}.

$\rhd$ Estimate \eqref{eq:PoincarLemma} is nothing but an accurate variant to the (sometimes called) ``Poincar\'e's Lemma" which is attributed to Mehler 1866 in \cite{CCLLV}, and has also been considered by many other authors, among them Poincar\'e, Borel, L\'evy, Sznitman, Diaconis-Freedman. 
 
$\rhd$  Estimate \eqref{eq:functLLN} in the case (iii) seems to be really new. It is a consequence of Theorem~\ref{theo:EquivKacChaos} together
with the following convergence result.

\begin{theo} [\cite{HaurayMischler,KleberSphere}] \label{theo:KacEntropRate}
For any $f \in \PPP_6(\R^d)$ with finite Fisher information, there exists a sequence $(F^N)$, where $F^N  := [f^{\otimes N}]_{\BSN}  \in \PPP(\BSN)$  is the ``conditioned to the Boltzmann's spheres of the product measures $f^{\otimes N}$", such that  
\bean
&&\bullet\quad W_1(F^N_\ell,f^{\otimes\ell}) \le {C_\ell \over N^{1/2}}, \,\,  \quad 1 \le \ell \le N,
\\
&&\bullet\quad |H(F^N|\gamma^N) - H(f|\gamma)|  \le {C \over N^{1/2}},
\\
&&\bullet\quad   I (F^N|\gamma^N)  \le C. 
\eean

\end{theo}

The proof just follows the one of a similar result (without rate) obtained in \cite{Kac1956,CCLLV} 
for the Kac's spheres framework and uses an accurate version (in $L^\infty$ norm) of the ``local central limit theorem''  firstly established by 
Berry and Esseen.

\subsection{Kac's, Boltmmann's and Fisher's chaos on Kac's and Boltzmann's spheres} 
In a last step, we establish some links between the different types of chaos on Kac's and Boltzmann's spheres.

\begin{theo}  [\cite{HaurayMischler,KleberSphere}]  \label{theo:KacEntropSphere}
For any sequence $G^N \in \Ps(E^N)$ and function $f \in \PPP(E)$, which fulfill some additional convenient 
moments and Fisher information bounds, there hold:
\[
|H(G^N|\gamma^N) - H(f|\gamma)| \le C \, I(G^N|\gamma^N) \, \Omega_j(G^N;f) + {C \over N^\alpha}, 
\]
\[
H(G^N_k|f^{\otimes k})  \le C_k \, I(G^N|\gamma^N) \, \Omega_j(G^N;f) + {C \over N^\alpha}. 
\]
\end{theo}

The proof of Theorem~\ref{theo:KacEntropSphere} is similar to the one of Theorem~\ref{theo:FBKchaos}
where however the HWI inequality is replaced by the following variant proved by Lott and Villani
$$
|H(G^N|\gamma^N) - H(F^N|\gamma^N)| \le C \, (I(G^N|\gamma^N)  + I(F^N|\gamma^N)) \, W_2(F^N,G^N) , 
$$
with $F^N := $  $[f^{\otimes N}]$. On the other hand, we use a result by 
Carlen-Lieb-Loss and Barthe-Cordero-Erausquin-Maurey in order to estimate $I(f|\gamma) \le \liminf I(G^N_1|\gamma^N_1) \le 
\liminf 2 I(G^N|\gamma^N)$ and next we use the bound on 
$|H(F^N|\gamma^N) - H(f|\gamma)|$ established in Theorem~\ref{theo:KacEntropRate}. 

\bigskip
\section{Conclusion and open problems}

 \subsection{Final step in the proof of Theorems~\ref{Th:MM-MM} and \ref{Th:MM-HS}}
We briefly explain how to deduce  Theorems~\ref{Th:MM-MM} and \ref{Th:MM-HS} and their consequences.

\medskip\noindent
{\bf Kac's Problem 1. }  By gathering Theorem~\ref{theo:Kac3} with the first estimate in Theorem~\ref{theo:KacEntropRate} (in order
to bound the term $T_3$ in the proof of Theorem~\ref{theo:Kac3}), 
with the last estimate with $j=N$ and $\ell = 2$ in Theorem~\ref{theo:EquivKacChaos} (in order to reformulate the LHS term in Theorem~\ref{theo:Kac3}) and some interpolation arguments (in order to estimate $W_1$ by $\| \cdot \|_{\FF'}$)
we obtain the first estimate in  Theorems~\ref{Th:MM-MM} and \ref{Th:MM-HS}. 
  
 \medskip\noindent
{\bf Kac's Problem 2. }  On the one hand, we write 
\bean
W_1(G^N_t,\gamma^N)
& \le& W_1(G^N_t, f_t^N) + W_1(f_t,\gamma) + W_1(\gamma^{\otimes N},\gamma^N)
\\
&\le& \Theta'(N) + C \, e^{-\lambda t} + {C \over \sqrt{N}}
\eean
as a consequence of Theorem~\ref{Th:MM-MM} or \ref{Th:MM-HS}, of \eqref{eq:ftogamma} and of Poincar\'e's lemma~\eqref{eq:PoincarLemma}. 
That gives a  good estimate for $N <\!< t$. 
We then optimize that last estimate together with  \eqref{eq:hNto1} which is a good estimate for $N>\!>t$ and we conclude to \eqref{eq:MM-GNtogammaN} and \eqref{eq:HS-GNtogammaN}. 

 \medskip\noindent
{\bf Kac's Problem 3. } For the {\bf M} model,  one can adapt Villani's result \cite{Villani:FisherInfoBoltzmann} for the Boltzmann equation to the Kac-Boltzmann equation 
and obtain (\cite{MaxHau,MMinvent}) the uniform estimate on the Fisher information
\beqn\label{eq:IFNtleIFN0}
\sup_{t \ge 0} I(G^N_t|\gamma^N) \le I(G^N_0|\gamma^N).
\eeqn
We then conclude thanks to Theorem~\ref{theo:KacEntropSphere}. 

For the {\bf HS} model, the proof is somewhat simpler, and in the same time less accurate.  We classically write 
 $$
 H(G^N_t|\gamma^N) + \int_0^t D(G^N_s) \, ds = H(G^N_0|\gamma^N)
 $$
 and
 $$
 H(f_t|\gamma) + \int_0^t D(f_s) \, ds = H(f_0|\gamma), 
 $$
 for some convenient dissipation of entropy terms $D(G^N)$ and $D(f)$. For the ``conditioned to the Boltzmann's spheres" $G^N_0$ initial datum 
 associated to $f_0$ we obtain by using some weak lsc property for the  entropy and dissipation of entropy functional the following series of inequalities
 \bean 
  H(f_t|\gamma) + \int_0^t D(f_s) \, ds
  &\le& \liminf_{N \to \infty} \Bigl\{ H(G^N_t|\gamma^N) + \int_0^t D(G^N_s) \, ds  \Bigr\}
 \\
&=& \liminf_{N \to \infty} H(G^N_0|\gamma^N) = H(f_0|\gamma)
\\
&=& H(f_t|\gamma) + \int_0^t D(f_s) \, ds,
\eean
 which in turns implies $H(G^N_t|\gamma^N) \to H(f_t|\gamma)$ as $N \to \infty$, for any $t \ge 0$.

 \subsection{Discussion} We just want to emphasize that the proof of the main Theorems~\ref{Th:MM-MM} \& \ref{Th:MM-HS}
 use tools from many different domains of mathematics, namely coming from  
 
 \smallskip
$\bullet$ functional analysis in increasing and infinite dimension,  and more specifically from theoretical statistics (functional LLN of Glivenko-Cantelli, 
 local central limit theorem of Berry-Esseen, mixtures according to De Finetti, Hewitt \& Savage), information theory (entropy of Boltzmann and Fisher information), optimal transport; 

 \smallskip
$\bullet$ probability theory (Markov semigroup); 
 
\smallskip$\bullet$ PDE (accurate stability estimates on the Boltzmann equation); 
 
 \smallskip$\bullet$ geometry (on the Kac's and Boltzmann's spheres).

 \subsection{Open problem:} 
 \

\smallskip
(1) Improve the rate of convergence for the Hard spheres model by improving the estimate for a weak distance between two solutions of the Boltzmann equation established in \cite{Fo-Mo}. 

\smallskip
(2) Obtain similar result for a true Hard potential (with rate) or even for the true soft potential (without rate to begin with). 
A first generalization of our method was obtained for the Landau equation (for Maxwell molecules) by Carrapatoso in \cite{KleberLandau}. 

\smallskip 
 (3) Consider singular models and generalize the propagation of chaos  for the vortex model obtained in  \cite{FHM}. 
 
\smallskip
 (4) Generalize the coupling technics used for  the asymmetric variant of  the Boltzmann-Kac system in \cite{FM}  to the Boltzmann-Kac system considered in these notes.

\bigskip
\bigskip
\bigskip
\bigskip

\bibliographystyle{acm}

  \signsm  

\end{document}